\documentclass[preprint,12pt,times]{elsarticle}

\usepackage{setspace}
\usepackage{geometry}
\usepackage{url}
\usepackage{graphicx,graphics}
\usepackage{subfigure}
\usepackage{epsfig}
\usepackage{psfrag}
\usepackage{hyperref}
\usepackage{amsmath,amssymb,amsfonts,amsthm}
\usepackage{mathrsfs}
\usepackage{color}
\usepackage{float}
\usepackage{natbib}
\usepackage{tabularx}
\usepackage{multirow}
\usepackage{array}

\theoremstyle{remark}

\newcommand{\Eref}[1]{Equation (\ref{#1})}
\newcommand{\fref}[1]{Figure (\ref{#1})}

\newcommand{\frefs}[1]{Figures~(\ref{#1})}

\newcommand{\xx}{\mathbf{x}}
\newcommand{\KK}{\mathbf{K}}
\newcommand{\bfm}{\mathbf{M}}
\newcommand{\bn}{\mathbf{N}}
\newcommand{\DD}{\mathbf{D_b}}

\newcommand{\bveps}{\boldsymbol{\varepsilon}}
\newcommand{\uu}{\mathbf{u}}
\newcommand{\BB}{\mathbf{B}}

\newcommand{\pp}{\mathbf{p}}
\newcommand{\bss}{\mathbf{s}}
\newcommand{\bdb}{\mathbf{b}}

\newcommand{\bds}{\boldsymbol{\delta}}

\newcommand{\rmd}{\mathrm{d}}

\newcommand\undermat[2]{%
  \makebox[0pt][l]{$\smash{\underbrace{\phantom{%
    \begin{matrix}#2\end{matrix}}}_{\text{$#1$}}}$}#2}

\journal{European Journal of Mechanics - A/Solids}

\begin{document}

\begin{frontmatter}

\title{Analysis of functionally graded material plates using triangular elements with cell-based smoothed discrete shear gap method}       

\author[unsw]{S Natarajan \corref{cor1}\fnref{fn1}}
\author[portugal,dubai]{AJM Ferreira}
\author[cardiff]{S Bordas}
\author[italy,dubai]{E Carrera}
\author[italy]{M Cinefra}
\author[egypt,dubai]{AM Zenkour}

\cortext[cor1]{Corresponding author}

\address[unsw]{School of Civil \& Environmental Engineering, The University of New South Wales, Sydney, Australia.}
\address[portugal]{Faculdade de Engenharia da Universidade do Porto, Porto, Portugal.}
\address[cardiff]{Institute of Mechanics and Advanced Materials, Cardiff School of Engineering, Cardiff University, Wales, UK.}
\address[italy]{Department of Aeronautics and Aerospace Engineering, Politecnico di Torino, Italy.}
\address[egypt]{Department of Mathematics, Faculty of Science, Kafrelsheikh University, Kafr El-Sheikh 33516, Egypt.}
\address[dubai]{Department of Mathematics, Faculty of Science, King Abdulaziz University, P.O. Box 80203, Jeddah 21589, Saudi Arabia.}
\fntext[fn1]{\url sundararajan.natarajan@gmail.com}

\begin{abstract}
In this paper, a cell based smoothed finite element method with discrete shear gap technique is employed to study the static bending, free vibration, mechanical and thermal buckling behaviour of functionally graded material (FGM) plates. The plate kinematics is based on the first order shear deformation theory and the shear locking is suppressed by a discrete shear gap method. The shear correction factors are evaluated by employing the energy equivalence principle. The material property is assumed to be temperature dependent and graded only in the thickness direction. The effective properties are computed by using the Mori-Tanaka homogenization method. The accuracy of the present formulation is validated against available solutions. A systematic parametric study is carried out to examine the influence the gradient index, the plate aspect ratio, skewness of the plate and the boundary conditions on the global response of the FGM plates. The effect of a centrally located circular cutout on the global response is also studied.
\end{abstract}


\end{frontmatter}


\section{Introduction}

With the rapid advancement of engineering,  there is an increasing demand for new materials which suits the harsh working environment without loosing it's mechanical, thermal or electrical properties. Engineered materials such as the composite materials are used due to their excellent strength-to and stiffness-to-weight ratios and their possibility of tailoring the properties in optimizing their structural response. But due to the abrupt change in material properties from matrix to fibre and between the layers, these materials suffer from pre-mature failure or by the decay in the stiffness characteristics because of delaminations and chemically unstable matrix and lamina adhesives. On the contrary, another class of materials, called, the Functionally Graded Materials (FGM) are made up of mixture of ceramics and metals and are characterized by \emph{smooth and continuous transition} in properties from one surface to another~\cite{koizumi1993}. As a result, FGMs are preferred over the laminated composites for structural integrity. The FGMs combine the best properties of the ceramics and the metals and this has attracted the researchers to study the characteristics of such structures. 

\paragraph{Background}
The tunable thermo-mechanical property of the FGM has attracted researchers to study the static and the dynamic behaviour of structures made of FGM under mechanical~\cite{Zenkour2006,Zenkour2007,reddy2000,singhatprakash2011} and thermal loading~\cite{natarajanbaiz2011, praveenreddy1998, dailiu2005, ganapathiprakash2006a, janghorbanazare2011, zenkourmashat2010, zhaolee2009}. Praveen \textit{et al.,}~\cite{praveenreddy1998} and Reddy \textit{et al.,}~\cite{reddychin2007} studied the thermo-elastic response of ceramic-metal plates using first order shear deformation theory coupled with the 3D heat conduction equation. Their study concluded that the structures made up of FGM with ceramic rich side exposed to elevated temperatures are susceptible to buckling due to the through thickness temperature variation. The buckling of skewed FGM plates under mechanical and thermal loads were studied in ~\cite{ganapathiprakash2006a,ganapathiprakash2006} employing the first order shear deformation theory and by using the shear flexible quadrilateral element. Efforts has also been made to study the mechanical behaviour of FGM plates with geometrical imperfection \cite{shariateslami2006}. Saji \textit{et al.,} \cite{sajivarughese2008} has studied thermal buckling of FGM plates with material properties dependent on both the composition and temperature. They found that Critical buckling temperatures are decreased when material properties are considered to be a function of temperature as compared to the results obtained where material properties are assumed to be independent of temperature. Ganapati~\textit{at al.,}~\cite{ganapathiprakash2006a} has studied the buckling of FGM skewed plate under thermal loading. Efforts has also been made to study the mechanical behaviour of FGM plates with geometrical imperfection \cite{shariateslami2006}. More recently, refined models have been adopted to study the characteristics of FGM structures~\cite{carrerabrischetto2008,carrrerabrischetto2011,cinefracarrera2012}.\\

Existing approaches in the literature to study plate and shell structures made up of FGMs uses finite element method (FEM) based on Lagrange basis functions~\cite{ganapathiprakash2006}, meshfree methods~\cite{ferreirabatra2006,qianbatra2004} and recently Valizadeh \textit{et al.,}~\cite{valizadehnatarajan2013} used non-uniform rational B-splines based FEM to study the static and dynamic characteristics of FGM plates in thermal environment. Even with these different approaches, the plate elements suffer from shear locking phenomenon and different techniques were proposed to alleviate the shear locking phenomenon.  Another set of methods have emerged to address the shear locking in the FEM. By incorporating the strain smoothing technique into the finite element method (FEM), Liu~\textit{et al.,}~\cite{liudai2007} have formulated a series of smoothed finite element methods (SFEM), named as cell-based SFEM (CS-FEM)~\cite{nguyenbordas2008,bordasnatarajan2010}, node-based SFEM~\cite{liunguyen2009}, edge-based SFEM~\cite{liunguyen2009a}, face-based SFEM~\cite{thoiliu2009} and $\alpha$-FEM~\cite{liunguyen2008}. And recently, edge based imbricate finite element method (EI-FEM) was proposed in~\cite{cazesmeschke2012} that shares common features with the ES-FEM. As the SFEM can be recast within a Hellinger-Reissner variational principle, suitable choices of the assumed strain/gradient space provides stable solutions. Depending on the number and geometry of the subcells used, a spectrum of methods exhibiting a spectrum of properties is obtained. Interested are referred to the literature~\cite{liudai2007,nguyenbordas2008} and references therein. Nguyen-Xuan~\textit{et al.,}~\cite{nguyenrabczuk2008} employed CS-FEM for Mindlin-Reissner plates. The curvature at each point is obtained by a non-local approximation via a smoothing function. From the numerical studies presented, it was concluded that the CS-FEM technique is robust, computationally inexpensive, free of locking and importantly insensitive to mesh distortions. The SFEM was extended to various problems such as shells~\cite{nguyenrabczuk2008a}, heat transfer~\cite{wuliu2010}, fracture mechanics~\cite{nguyenliu2013} and structural acoustics~\cite{hecheng2011} among others. In~\cite{bordasnatarajan2011}, CS-FEM has been combined with the extended FEM to address problems involving discontinuities. The above list is no way comprehensive and interested readers are referred to the literature and references therein and a recent review paper by Jha and Kant~\cite{jhakant2013} on FGM plates. 

\paragraph{Approach} In this paper, we study the static and the dynamic characteristics of FGM plates by using a cell-based smoothed finite element method with discrete shear gap technique. Three-noded triangular element is employed in this study. The effect of different parameters viz., the material gradient index, the plate aspect ratio, the plate slenderness ratio and the boundary conditon on the global response of FGM plates are numerically studied. The effect of centrally located circular cutout is also studied. The present work focusses on the computational aspects of the governing equations, hence, the attention has been restricted to Reissner-Mindlin plate theories. It is noted that the extension to higher order theories is possible.

\paragraph{Outline} The paper is organized as follows, the next section will give an introduction to FGM and a brief overview of Reissner-Mindlin plate theory. Section \ref{csdsg3} presents an overview of the cell-based smoothed finite element method with discrete shear gap technique. The efficiency of the present formulation, numerical results and parametric studies are presented in Section \ref{numresu}, followed by concluding remarks in the last section.

\section{Theoretical Background}
\label{theory}
\subsection{Reissner-Mindlin plate theory}
\begin{figure}[htpb]
\centering
\subfigure[]{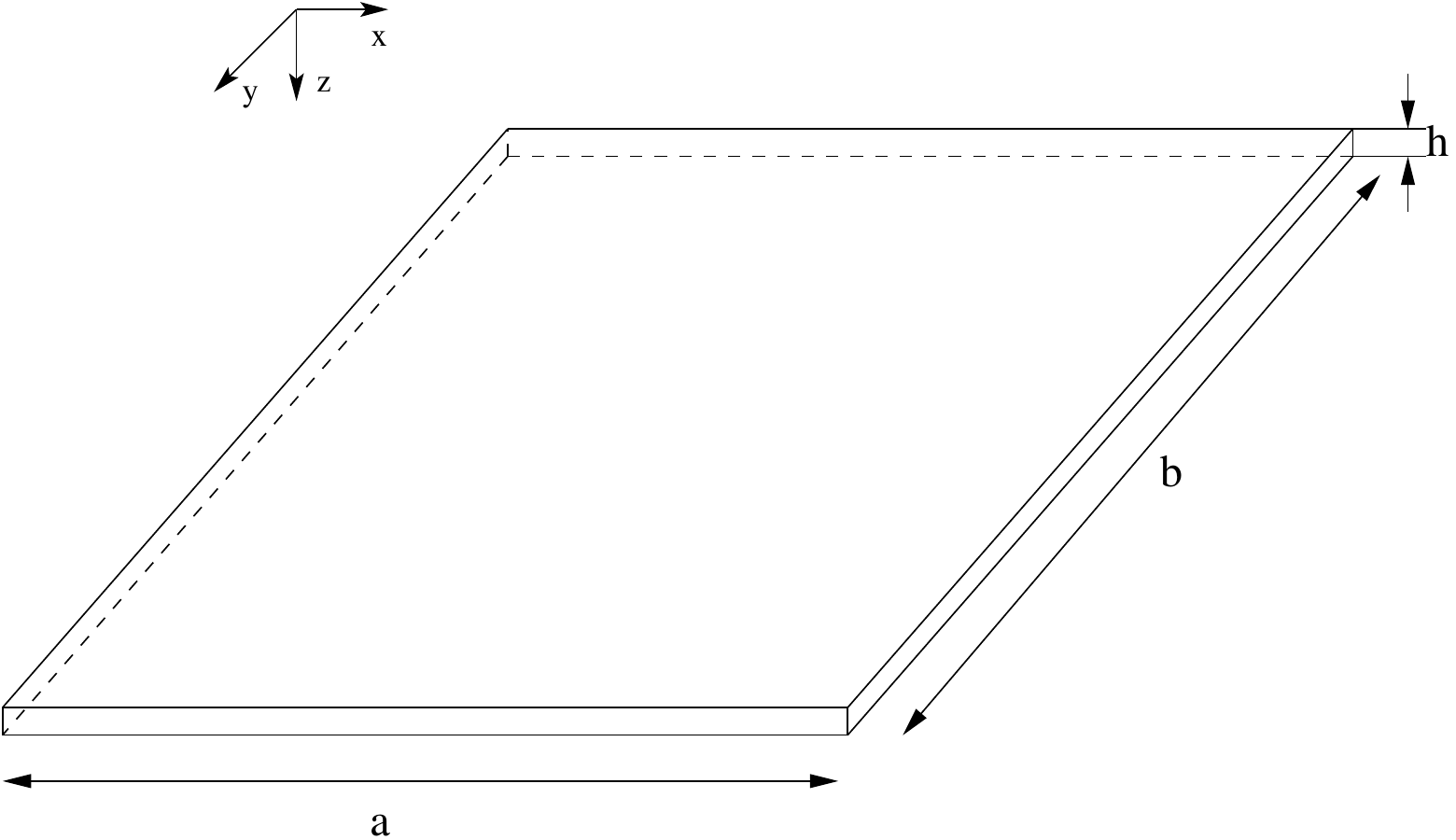}
\subfigure[]{\input{./Figures/skew.pstex_t}}
\caption{(a) coordinate system of a rectangular FGM plate, (b) Coordinate system of a skew plate}
\label{fig:platefig}
\end{figure}

The Reissner-Mindlin plate theory, also known as the first order shear deformation theory (FSDT) takes into account the shear deformation through the thickness. Using the Mindlin formulation, the displacements $u,v,w$ at a point $(x,y,z)$ in the plate (see \fref{fig:platefig}) from the medium surface are expressed as functions of the mid-plane displacements $u_o,v_o,w_o$ and independent rotations $\theta_x,\theta_y$ of the normal in $yz$ and $xz$ planes, respectively, as:

\begin{eqnarray}
u(x,y,z,t) &=& u_o(x,y,t) + z \theta_x(x,y,t) \nonumber \\
v(x,y,z,t) &=& v_o(x,y,t) + z \theta_y(x,y,t) \nonumber \\
w(x,y,z,t) &=& w_o(x,y,t) 
\label{eqn:displacements}
\end{eqnarray}
where $t$ is the time. The strains in terms of mid-plane deformation can be written as:

\begin{equation}
\bveps  = \left\{ \begin{array}{c} \bveps_p \\ 0 \end{array} \right \}  + \left\{ \begin{array}{c} z \bveps_b \\ \bveps_s \end{array} \right\} 
\label{eqn:strain1}
\end{equation}

The midplane strains $\bveps_p$, the bending strain $\bveps_b$ and the shear strain $\varepsilon_s$ in \Eref{eqn:strain1} are written as:

\begin{eqnarray}
\renewcommand{\arraystretch}{1.2}
\bveps_p = \left\{ \begin{array}{c} u_{o,x} \\ v_{o,y} \\ u_{o,y}+v_{o,x} \end{array} \right\}, \hspace{1cm}
\renewcommand{\arraystretch}{1.2}
\bveps_b = \left\{ \begin{array}{c} \theta_{x,x} \\ \theta_{y,y} \\ \theta_{x,y}+\theta_{y,x} \end{array} \right\}, \nonumber \\
\renewcommand{\arraystretch}{1.2}
\bveps_s = \left\{ \begin{array}{c} \theta _x + w_{o,x} \\ \theta _y + w_{o,y} \end{array} \right\}. \hspace{1cm}
\renewcommand{\arraystretch}{1.2}
\end{eqnarray}
where the subscript `comma' represents the partial derivative with respect to the spatial coordinate succeeding it. The membrane stress resultants $\bn$ and the bending stress resultants $\bfm$ can be related to the membrane strains, $\bveps_p$ and bending strains $\bveps_b$ through the following constitutive relations:

\begin{eqnarray}
\bn &=& \left\{ \begin{array}{c} N_{xx} \\ N_{yy} \\ N_{xy} \end{array} \right\} = \mathbf{A} \bveps_p + \BB \bveps_b - \bn^{\rm th}\nonumber \\
\bfm &=& \left\{ \begin{array}{c} M_{xx} \\ M_{yy} \\ M_{xy} \end{array} \right\} = \BB \bveps_p + \DD \bveps_b - \bfm^{\rm th}
\end{eqnarray}
where the matrices $\mathbf{A} = A_{ij}, \BB= B_{ij}$ and $\DD = D_{ij}; (i,j=1,2,6)$ are the extensional, bending-extensional coupling and bending stiffness coefficients and are defined as:

\begin{equation}
\left\{ A_{ij}, ~B_{ij}, ~ D_{ij} \right\} = \int_{-h/2}^{h/2} \overline{Q}_{ij} \left\{1,~z,~z^2 \right\}~dz
\end{equation}
Similarly, the transverse shear force $Q = \{Q_{xz},Q_{yz}\}$ is related to the transverse shear strains $\varepsilon_s$ through the following equation:

\begin{equation}
Q_{ij} = E_{ij} \varepsilon_s
\end{equation}

where $E_{ij} = \int_{-h/2}^{h/2} \overline{Q}_{ij} \upsilon_i \upsilon_j~dz;~ (i,j=4,5)$ are the transverse shear stiffness coefficients, $\upsilon_i, \upsilon_j$ is the transverse shear coefficient for non-uniform shear strain distribution through the plate thickness. The stiffness coefficients $\overline{Q}_{ij}$ are defined as:

\begin{eqnarray}
\overline{Q}_{11} = \overline{Q}_{22} = {E(z) \over 1-\nu^2}; \hspace{1cm} \overline{Q}_{12} = {\nu E(z) \over 1-\nu^2}; \hspace{1cm} \overline{Q}_{16} = \overline{Q}_{26} = 0 \nonumber \\
\overline{Q}_{44} = \overline{Q}_{55} = \overline{Q}_{66} = {E(z) \over 2(1+\nu) }
\end{eqnarray}

where the modulus of elasticity $E(z)$ and Poisson's ratio $\nu$ are given by \Eref{eqn:young}. The thermal stress resultant $\bn^{\rm th}$ and the moment resultant $\bfm^{\rm th}$ are:

\begin{eqnarray}
\bn^{\rm th}&=& \left\{ \begin{array}{c} N^{\rm th}_{xx} \\ N^{\rm th}_{yy} \\ N^{\rm th}_{xy} \end{array} \right\} = \int\limits_{-h/2}^{h/2} \overline{Q}_{ij} \alpha(z,T) \left\{ \begin{array}{c} 1 \\ 1 \\ 0 \end{array} \right\} \Delta T(z)~ \rmd z \nonumber \\
\bfm^{\rm th} &=& \left\{ \begin{array}{c} M^{\rm th}_{xx} \\ M^{\rm th}_{yy} \\ M^{\rm th}_{xy} \end{array} \right\} = \int\limits_{-h/2}^{h/2} \overline{Q}_{ij} \alpha(z,T) \left\{ \begin{array}{c} 1 \\ 1 \\ 0 \end{array} \right\} \Delta T(z)~ z ~\rmd z \nonumber \\ 
\end{eqnarray}

where the thermal coefficient of expansion $\alpha(z,T)$ is given by \Eref{eqn:thermalcondalpha} and $\Delta T(z) = T(z)-T_o$ is the temperature rise from the reference temperature and $T_o$ is the temperature at which there are no thermal strains. The strain energy function $U$ is given by:

\begin{equation}
\begin{split}
U(\boldsymbol{\delta}) = {1 \over 2} \int_{\Omega} \left\{ \bveps_p^{\textup{T}} \mathbf{A} \bveps_p + \bveps_p^{\textup{T}} \mathbf{B} \bveps_b + 
\bveps_b^{\textup{T}} \mathbf{B} \bveps_p + \bveps_b^{\textup{T}} \mathbf{D} \bveps_b +  \bveps_s^{\textup{T}} \mathbf{E} \bveps_s - \bveps_b^{\rm T} \bn^{\rm th} - \bveps_b^{\rm T} \bfm^{\rm th} \right\}~ \mathrm{d} \Omega
\end{split}
\label{eqn:potential}
\end{equation}

where $\boldsymbol{\delta} = \{u,v,w,\theta_x,\theta_y\}$ is the vector of the degree of freedom associated to the displacement field in a finite element discretization. Following the procedure given in~\cite{Rajasekaran1973}, the strain energy function $U$ given in~\Eref{eqn:potential} can be rewritten as:

\begin{equation}
U(\boldsymbol{\delta}) = {1 \over 2}  \boldsymbol{\delta}^{\textup{T}} \mathbf{K}  \boldsymbol{\delta}
\label{eqn:poten}
\end{equation}

where $\mathbf{K}$ is the linear stiffness matrix. The kinetic energy of the plate is given by:

\begin{equation}
T(\boldsymbol{\delta}) = {1 \over 2} \int_{\Omega} \left\{p (\dot{u}_o^2 + \dot{v}_o^2 + \dot{w}_o^2) + I(\dot{\theta}_x^2 + \dot{\theta}_y^2) \right\}~\mathrm{d} \Omega
\label{eqn:kinetic}
\end{equation}

where $p = \int_{-h/2}^{h/2} \rho(z)~dz, ~ I = \int_{-h/2}^{h/2} z^2 \rho(z)~dz$ and $\rho(z)$ is the mass density that varies through the thickness of the plate. When the plate is subjected to a temperature field, this in turn results in in-plane stress resultants, $\bn^{\rm th}$. The external work due to the in-plane stress resultants developed in the plate under a thermal load is given by:

\begin{equation}
\begin{split}
V(\boldsymbol{\delta}) = \int\limits_\Omega \left\{ \frac{1}{2} \left[ N_{xx}^{\rm th} w_{,x}^2 + N_{yy}^{\rm th} w_{,y}^2 + 2 N_{xy}^{\rm th}w_{,x}w_{,y}\right] + \right. \\ \left.
\frac{h^2}{24} \left[ N_{xx}^{\rm th} \left( \theta_{x,x}^2 + \theta_{y,x}^2 \right) + N_{yy}^2 \left( \theta_{x,y}^2 + \theta_{y,y}^2 \right) + 2 N_{xy}^{\rm th} \left( \theta_{x,x}\theta_{x,y} + \theta_{y,x}\theta_{y,y} \right) \right] \right\}~ d\Omega
\end{split}
\label{eqn:V}
\end{equation}

Substituting \Eref{eqn:potential} - (\ref{eqn:V}) in Lagrange's equation of motion, one obtains the following finite element equations:

{\bf Static bending}: 
\begin{equation}
\KK \boldsymbol{\delta} = \mathbf{F}
\label{eqn:staticgovern}
\end{equation}

{\bf Free vibration}: 
\begin{equation}
\bfm \ddot{\boldsymbol{\delta}} + \left( \KK + \KK_G \right) \boldsymbol{\delta} = \mathbf{0}
\label{eqn:freevibgovern}
\end{equation}

{\bf Buckling analysis}:

\paragraph*{Mechanical Buckling\footnote{Prebuckling deformations are assumed to be zero or negligible in the analysis (including those coming from in-plane and out-of-plane coupling related to FGM and temperature variation through the thickness of the plate).}}
\begin{equation}
\left( \KK + \lambda_M \KK_G \right) \boldsymbol{\delta} = \mathbf{0}
\end{equation}

\paragraph*{Thermal Buckling}
\begin{equation}
\left( \KK + \lambda_T \KK_G \right) \boldsymbol{\delta} = \mathbf{0}
\end{equation}

where $\lambda_M$  includes  the critical  value applied in-plane mechanical loading and $\lambda_T$ is the critical temperature difference and $\KK$, $\KK_G$ are the linear stiffness and geometric stiffness matrices, respectively. The critical temperature difference is computed using a standard eigenvalue algorithm.

\subsection{Functionally graded material}
A rectangular plate made of a mixture of ceramic and metal is considered with the coordinates $x,y$ along the in-plane directions and $z$ along the thickness direction (see \fref{fig:platefig}). The material on the top surface $(z=h/2)$ of the plate is ceramic rich and is graded to metal at the bottom surface of the plate $(z=-h/2)$ by a power law distribution. The effective properties of the FGM plate can be computed by using the rule of mixtures or by employing the Mori-Tanaka homogenization scheme. Let $V_i (i=c,m)$ be the volume fraction of the phase material. The subscripts $c$ and $m$ refer to ceramic and metal phases, respectively. The volume fraction of ceramic and metal phases are related by $V_c + V_m = 1$ and $V_c$ is expressed as:

\begin{equation}
V_c(z) = \left( \frac{2z+h}{2h} \right)^n
\end{equation}

where $n$ is the volume fraction exponent $(n \geq 0)$, also known as the gradient index. The variation of the composition of ceramic and metal is linear for $n=$1, the value of $n=$ 0 represents a fully ceramic plate and any other value of $n$ yields a composite material with a smooth transition from ceramic to metal.

\paragraph*{Mori-Tanaka homogenization method}
Based on the Mori-Tanaka homogenization method, the effective Young's modulus and Poisson's ratio are computed from the effective bulk modulus $K$ and the effective shear modulus $G$ as~\cite{Sundararajan2005}

\begin{equation}
\frac{K_{\rm eff}-K_m}{K_c-K_m} = \frac{V_c}{1 + V_m \frac{3(K_c-K_m)}{3K_m+4G_m}}, \hspace{0.75cm}
\frac{G_{\rm eff}-G_m}{G_c-G_m} = \frac{V_c}{1 + V_m \frac{(G_c-G_m)}{(G_m+f_1)} }
\end{equation}
where

\begin{equation}
f_1 = \frac{G_m (9K_m+8G_m)}{6(K_m+2G_m)}
\end{equation}

The effective Young's modulus $E_{\rm eff}$ and Poisson's ratio $\nu_{\rm eff}$ can be computed from the following relations:

\begin{equation}
E_{\rm eff} = \frac{9 K_{\rm eff} G_{\rm eff}}{3K_{\rm eff} + G_{\rm eff}}, \hspace{1cm} \nu_{\rm eff} = \frac{3K_{\rm eff} - 2G_{\rm eff}}{2(3K_{\rm eff} + G_{\rm eff})}
\label{eqn:young}
\end{equation}

The effective mass density $\rho_{\rm eff}$ is computed using the rule of mixtures $(\rho_{eff} = \rho_cV_c + \rho_mV_m)$. The effective heat conductivity $\kappa_{\rm eff}$ and the coefficient of thermal expansion $\alpha_{\rm eff}$ is given by:

\begin{eqnarray}
\frac{\kappa_{\rm eff} - \kappa_m}{\kappa_c - \kappa_m} = \frac{V_c}{1 + V_m \frac{(\kappa_c - \kappa_m)}{3\kappa_m}} \nonumber \\
\frac{\alpha_{\rm eff} - \alpha_m}{\alpha_c - \alpha_m} = \frac{ \left( \frac{1}{K_{\rm eff}} - \frac{1}{K_m} \right)}{\left(\frac{1}{K_c} - \frac{1}{K_m} \right)}
\label{eqn:thermalcondalpha}
\end{eqnarray}

\paragraph*{Temperature dependent material property} The material properties that are temperature dependent are written as~\cite{Sundararajan2005}:
\begin{equation}
P = P_o(P_{-1}T^{-1} + 1 + P_1T + P_2T^2 + P_3T^3)
\end{equation}
where $P_o,P_{-1},P_1,P_2$ and $P_3$ are the coefficients of temperature $T$ and are unique to each constituent material phase.
\paragraph*{Temperature distribution through the thickness}  The temperature variation is assumed to occur in the thickness direction only and the temperature field is considered to be constant in the $xy$-plane. In such a case, the temperature distribution along the thickness can be obtained by solving a steady state heat transfer problem:

\begin{equation}
-{d \over dz} \left[ \kappa(z) {dT \over dz} \right] = 0, \hspace{0.5cm} T = T_c ~\textup{at}~ z = h/2;~~ T = T_m ~\textup{at} ~z = -h/2
\label{eqn:heat}
\end{equation}

The solution of \Eref{eqn:heat} is obtained by means of a polynomial series~\cite{wu2004} as

\begin{equation}
T(z) = T_m + (T_c - T_m) \eta(z,h)
\label{eqn:tempsolu}
\end{equation}

where,

\begin{equation}
\begin{split}
\eta(z,h) = {1 \over C} \left[ \left( {2z + h \over 2h} \right) - {\kappa_{cm} \over (n+1)\kappa_m} \left({2z + h \over 2h} \right)^{n+1} + \right. \\ 
\left. {\kappa_{cm} ^2 \over (2n+1)\kappa_m ^2 } \left({2z + h \over 2h} \right)^{2n+1}
-{\kappa_{cm} ^3 \over (3n+1)\kappa_m ^3 } \left({2z + h \over 2h} \right)^{3n+1} \right. \\ + 
\left. {\kappa_{cm} ^4 \over (4n+1)\kappa_m^4 } \left({2z + h \over 2h} \right)^{4n+1} 
- {\kappa_{cm} ^5 \over (5n+1)\kappa_m ^5 } \left({2z + h \over 2h} \right)^{5n+1} \right] ;
\end{split}
\label{eqn:heatconducres}
\end{equation}

\begin{equation}
\begin{split}
C = 1 - {\kappa_{cm} \over (n+1)\kappa_m} + {\kappa_{cm} ^2 \over (2 n+1)\kappa_m ^2} 
- {\kappa_{cm} ^3 \over (3n+1)\kappa_m ^3} \\ + {\kappa_{cm} ^4 \over (4n+1)\kappa_m ^4}
- {\kappa_{cm} ^5\over (5n+1)\kappa_m ^5}
\end{split}
\end{equation}

\section{Cell based smoothed finite element method with discrete shear gap technique}\label{csdsg3}
In this study, three-noded triangular element with five degrees of freedom (dofs) $\boldsymbol{\delta} = \{u,v,w,\theta_x,\theta_y\}$ is employed. The displacement is approximated by
\begin{equation}
\uu^h = \sum_IN_I \boldsymbol{\delta}_I
\end{equation}
where $\boldsymbol{\delta}_I$ are the nodal dofs and $N_I$ are the standard finite element shape functions given by
\begin{equation}
N = \left[ 1-\xi-\eta, \;\; \eta, \;\; \xi \right]
\end{equation}
\begin{figure}[htpb]
\centering
\scalebox{0.8}{\input{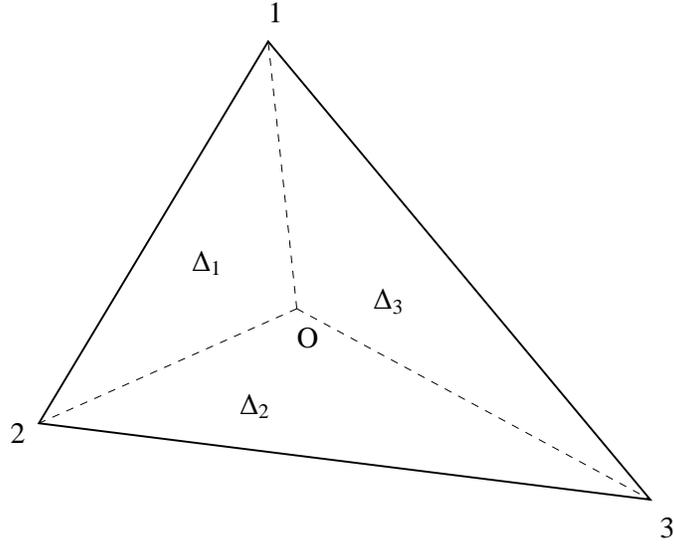}}
\caption{A triangular element is divided into three subtriangles. $\Delta_1, \Delta_2$ and $\Delta_3$ are the subtriangles created by connecting the central point $O$ with three field nodes.}
\label{fig:triEle}
\end{figure}
In the proposed approach, cell-based smoothed finite element method (CSFEM) is combined with stabilized discrete shear gap method (DSG) for three-noded triangular element, called as `cell-based discrete shear gap method (CS-DSG3).' The cell-based smoothing technique decreases the computational complexity, whilst DSG suppresses the shear locking phenomenon when the present formulation is applied to thin plates. Interested readers are referred to the literature and references therein for the description of cell-based smoothing technique~\cite{liudai2007,bordasnatarajan2010} and DSG method~\cite{bletzingerbischoff2000}. In the CS-DSG3, each triangular element is divided into three subtriangles. The displacement vector at the center node is assumed to be the simple average of the three displacement vectors of the three field nodes. In each subtriangle, the stabilized DSG3 is used to compute the strains and also to avoid the transverse shear locking. Then the strain smoothing technique on the whole triangular element is used to smooth the strains on the three subtriangles.  Consider a typical triangular element $\Omega_e$ as shown in \fref{fig:triEle}. This is first divided into three subtriangles $\Delta_1, \Delta_2$ and $\Delta_3$ such that $\Omega_e = \bigcup\limits_{i=1}^3 \Delta_i$. The coordinates of the center point $\xx_o = (x_o,y_o)$ is given by:
\begin{equation}
(x_o,y_o) = \frac{1}{3} (x_I, y_I)
\end{equation}
The displacement vector of the center point is assumed to be a simple average of the nodal displacements as
\begin{equation}
\boldsymbol{\delta}_{eO} = \frac{1}{3}\boldsymbol{\delta}_{eI}
\label{eqn:centerdefl}
\end{equation}

The constant membrane strains, the bending strains and the shear strains for subtriangle $\Delta_1$ is given by:
\begin{align}
\bveps_p &= \left[ \begin{array}{ccc} \pp_1^{\Delta_1} & \pp_2^{\Delta_1} & \pp_3^{\Delta_1} \end{array} \right] \left\{ \begin{array}{c} \bds_{eO} \\ \bds_{e1} \\ \bds_{e2} \end{array} \right\} \nonumber \\
\bveps_b &= \left[ \begin{array}{ccc} \bdb_1^{\Delta_1} & \bdb_2^{\Delta_1} & \bdb_3^{\Delta_1} \end{array} \right] \left\{ \begin{array}{c} \bds_{eO} \\ \bds_{e1} \\ \bds_{e2} \end{array} \right\} \nonumber \\
\bveps_s &= \left[ \begin{array}{ccc} \bss_1^{\Delta_1} & \bss_2^{\Delta_1} & \bss_3^{\Delta_1} \end{array} \right] \left\{ \begin{array}{c} \bds_{eO} \\ \bds_{e1} \\ \bds_{e2} \end{array} \right\}
\label{eqn:constrains}
\end{align}

\begin{figure}[htpb]
\centering
\scalebox{0.7}{\input{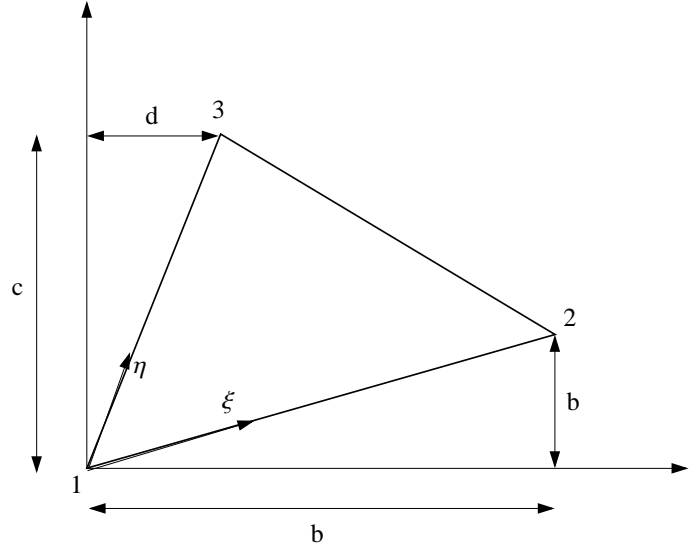}}
\caption{Three-noded triangular element and local coordinates in discrete shear gap method.}
\label{fig:dsg3}
\end{figure}

Upon substituting the expression for $\boldsymbol{\delta}_{eO}$ in \Eref{eqn:constrains}, we obtain:

\begin{align}
\bveps_p ^{\Delta_1}&= \left[ \begin{array}{ccc} \frac{1}{3}\pp_1^{\Delta_1} + \pp_2^{\Delta_1}& \frac{1}{3}\pp_1^{\Delta_1} + \pp_3^{\Delta_1} & \frac{1}{3} \pp_1^{\Delta_1} \end{array} \right] \left\{ \begin{array}{c} \bds_{e1} \\ \bds_{e2} \\ \bds_{e3} \end{array} \right\} = \BB_p^{\Delta_1} \bds_e \nonumber \\
\bveps_b^{\Delta_1} &= \left[ \begin{array}{ccc} \frac{1}{3}\bdb_1^{\Delta_1} + \bdb_2^{\Delta_1}& \frac{1}{3}\bdb_1^{\Delta_1} + \bdb_3^{\Delta_1} & \frac{1}{3} \bdb_1^{\Delta_1} \end{array} \right] \left\{ \begin{array}{c} \bds_{e1} \\ \bds_{e2} \\ \bds_{e3} \end{array} \right\} = \BB_b^{\Delta_1} \bds_e\nonumber \\
\bveps_s^{\Delta_1} &= \left[ \begin{array}{ccc} \frac{1}{3}\bss_1^{\Delta_1} + \bss_2^{\Delta_1}& \frac{1}{3}\bss_1^{\Delta_1} + \bss_3^{\Delta_1} & \frac{1}{3} \bss_1^{\Delta_1} \end{array} \right] \left\{ \begin{array}{c} \bds_{e1} \\ \bds_{e2} \\ \bds_{e3} \end{array} \right\} = \BB_s^{\Delta_1} \bds_e \nonumber \\
\label{eqn:constrainsSubtri}
\end{align}

where $\pp_i, (i=1,2,3)$, $\bdb_i, (i=1,2,3)$ and $\bss_i, (i=1,2,3)$ are given by:
\begin{align}
\BB_p &= \frac{1}{2A_e} \left[ \begin{array}{rrrrrrrrrrrrrrr} b-c & 0 & 0 & 0 & 0 & c & 0 & 0 & 0 & 0 & -b & 0 & 0 & 0 & 0 \\ 0 & d-a & 0 & 0 & 0 & 0 & -d & 0 & 0 & 0 & a & 0 & 0 & 0 & 0 \\\undermat{\pp_1}{ d-a & b-c & 0 & 0 & 0} & \undermat{\pp_2}{-d & c & 0 & 0 & 0} & \undermat{\pp_3}{a & -b & 0 & 0 & 0} \end{array} \right] \nonumber \\ \nonumber \\ \nonumber \\
\BB_b & = \frac{1}{2A_e} \left[ \begin{array}{rrrrrrrrrrrrrrr} 0 & 0 & 0 & b-c & 0 & 0 & 0 & 0 & c & 0 & 0 & 0 & 0 & -b & 0 \\ 0 & 0 & 0 & 0 & d-a & 0 & 0 & 0 & 0 & -d & 0 & 0 & 0 & 0 & a \\ \undermat{\bdb_1}{0 & 0 & 0 & d-a & b-c} & \undermat{\bdb_2}{0 & 0 & 0 & -d & c} & \undermat{\bdb_3}{0 & 0 & 0 & a & -b} \end{array} \right] \nonumber \\ \nonumber \\ \nonumber \\
\BB_s & = \frac{1}{2A_e} \left[ \begin{array}{rrrrrrrrrrrrrrr} 0 & 0 & b-c & A_e & 0 & 0 & 0 & c & ac/2 & bc/2 & 0 & 0 & -b & -bd/2 & -bc/2 \\ \undermat{\bss_1}{0 & 0 & d-a & 0 & A_e} & \undermat{\bss_2}{0 & 0 & -d & -ad/2 & -bd/2} & \undermat{\bss_3}{0 & 0 & a & ad/2 & ac/2} \end{array} \right] \\
\end{align}

where $a = x_2 - x_1; b = y_2 - y_1; c = y_3 - y_1$ and $d = x_3 - x_1$ (see \fref{fig:dsg3}), $A_e$ is the area of the triangular element and $\BB_s$ is altered shear strains~\cite{bletzingerbischoff2000}. The strain-displacement matrix for the other two triangles can be obtained by cyclic permutation.

Now applying the cell-based strain smoothing~\cite{bordasnatarajan2010}, the constant membrane strains, the bending strains and the shear strains are respectively employed to create a smoothed membrane strain $\overline{\bveps}_p$, smoothed bending strain $\overline{\bveps}_b$ and smoothed shear strain $\overline{\bveps}_s$on the triangular element $\Omega_e$ as:

\begin{align}
\overline{\bveps}_p &= \int\limits_{\Omega_e} \bveps_b \Phi_e(\xx)~\rmd \Omega = \sum\limits_{i=1}^3 \bveps_p^{\Delta_i} \int\limits_{\Delta_i} \Phi_e(\xx)~\rmd \Omega \nonumber \\
\overline{\bveps}_b &= \int\limits_{\Omega_e} \bveps_b \Phi_e(\xx)~\rmd \Omega = \sum\limits_{i=1}^3 \bveps_b^{\Delta_i} \int\limits_{\Delta_i} \Phi_e(\xx)~\rmd \Omega \nonumber \\
\overline{\bveps}_s &= \int\limits_{\Omega_e} \bveps_s \Phi_e(\xx)~\rmd \Omega = \sum\limits_{i=1}^3 \bveps_s^{\Delta_i} \int\limits_{\Delta_i} \Phi_e(\xx)~\rmd \Omega
\end{align}

where $\Phi_e(\xx)$ is a given smoothing function that satisfies. In this study, following constant smoothing function is used:
\begin{equation}
\Phi(\xx) = \left\{ \begin{array}{cc} 1/A_c & \xx \in \Omega_c \\ 0 & \xx \notin \Omega_c \end{array} \right.
\end{equation}
where $A_c$ s the area of the triangular element, the smoothed membrane strain, the smoothed bending strain and the smoothed shear strain is then given by
\begin{equation}
\left\{ \overline{\bveps}_p, \overline{\bveps}_b, \overline{\bveps}_s \right\} = \frac{ \sum\limits_{i=1}^3 A_{\Delta_i} \left\{ \bveps_p^{\Delta_i}, \bveps_b^{\Delta_i}, \bveps_s^{\Delta_i} \right \} }{A_e}
\end{equation}

The smoothed elemental stiffness matrix is given by
\begin{align}
\KK &= \int\limits_{\Omega_e} \overline{\BB}_p \mathbf{A} \overline{\BB}_p^{\rm T} + \overline{\BB}_p \BB \overline{\BB}_b^{\rm T} + \overline{\BB}_b \BB \overline{\BB}_p^{\rm T} + \overline{\BB}_b \mathbf{D} \overline{\BB}_b^{\rm T} + \overline{\BB}_s \mathbf{E} \overline{\BB}_s^{\rm T}~\rmd \Omega \nonumber \\
&= \left( \overline{\BB}_p \mathbf{A} \overline{\BB}_p^{\rm T} + \overline{\BB}_p \BB \overline{\BB}_b^{\rm T} + \overline{\BB}_b \BB \overline{\BB}_p^{\rm T} + \overline{\BB}_b \mathbf{D} \overline{\BB}_b^{\rm T} + \overline{\BB}_s \mathbf{E} \overline{\BB}_s^{\rm T} \right) A_e
\end{align}
where $\overline{\BB}_p, \overline{\BB}_b$ and $\overline{\BB}_s$ are the smoothed strain-displacement matrix. The mass matrix $\bfm$, is computed by following the conventional finite element procedure.
\section{Numerical examples} \label{numresu}
In this section, we present the static bending response, the linear free vibration and buckling analysis of FGM plates using cell based smoothed finite element method with discrete shear gap technique. The effect of various parameters, viz., material gradient index $n$, skewness of the plate $\psi$, the plate aspect ratio $a/b$, the plate thickness $a/h$ and boundary conditions on the global response is numerically studied. The top surface of the plate is ceramic rich and the bottom surface of the plate is metal rich. Here, the modified shear correction factor obtained based on energy equivalence principle as outlined in~\cite{singhaprakash2011} is used. The boundary conditions for simply supported and clamped cases are :
\noindent \emph{Simply supported boundary condition}: \\
\begin{equation}
u_o = w_o = \theta_y = 0 \hspace{0.2cm} ~\textup{on} \hspace{0.2cm}  x=0,a; \hspace{0.2cm}
v_o = w_o = \theta_x = 0 \hspace{0.2cm} ~\textup{on} \hspace{0.2cm}  y=0,b
\end{equation}

\noindent \emph{Clamped boundary condition}: \\
\begin{equation}
u_o = w_o = \theta_y = v_o = \theta_x = 0 \hspace{0.25cm} ~\textup{on} ~ x=0,a \hspace{0.2cm} \& \hspace{0.2cm}  y=0,b
\end{equation}

\paragraph*{Skew boundary transformation} For skew plates, the edges of the boundary elements may not be parallel to the global axes $(x,y,z)$. In order to specify the boundary conditions on skew edges, it is necessary to use the edge displacements $(u_o^\prime,v_o^\prime,w_o^\prime)$ etc, in a local coordinate system $(x^\prime,y^\prime,z^\prime)$ (see \fref{fig:platefig}). The element matices corresponding to the skew edges are transformed from global axes to local axes on which the boundary conditions can be conveniently specified. The relation between the global and the local degrees of freedom of a particular node is obtained by:
\begin{equation}
\boldsymbol{\delta} = \mathbf{L}_g \boldsymbol{\delta}^\prime
\end{equation}
where $\boldsymbol{\delta}$ and $\boldsymbol{\delta}^\prime$ are the generalized displacement vector in the global and the local coordinate system, respectively. The nodal transformation matrix for a node $I$ on the skew boundary is given by:
\begin{equation}
\mathbf{L}_g = \left[ \begin{array}{rrrrr} \cos\psi & \sin\psi & 0 & 0 & 0 \\ -\sin\psi & \cos\psi & 0 & 0 & 0 \\ 0 & 0 & 1 & 0 & 0 \\ 0 & 0 & 0 & \cos\psi & \sin\psi \\ 0 & 0 & 0 & -\sin\psi & \cos\psi \end{array} \right]
\end{equation}
where $\psi$ defines the skewness of the plate.

\subsection{Static Bending}

\begin{figure}[htpb]
\centering
\includegraphics[scale=0.7]{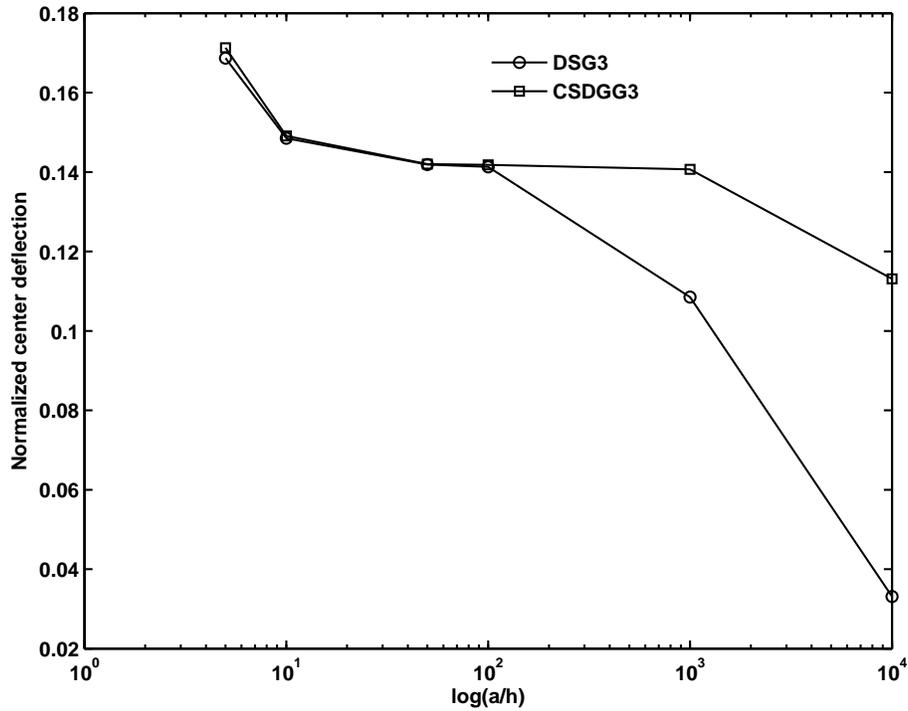}
\caption{The normalized center deflection as a function of normalized plate thickness for a simply supported square FGM plate subjected to a uniform load.}
\label{fig:plateshearlock}
\end{figure}

Let us consider a Al/ZrO$_2$ FGM square plate with length-to-thickness $a/h=$ 5, subjected to a uniform load with fully simply supported (SSSS) and fully clamped (CCCC) boundary conditions. The Young's modulus for ZrO$_2$ is $E_c=$ 151 GPa and for aluminum is $E_m=$ 70 GPa. Poisson's ratio is chosen as constant, $\nu=$ 0.3. Table \ref{tab:ssss-valid} compares the results from the present formulation with other approaches available in the literature~\cite{gilhooleybatra2007,leezhao2009,nguyen-xuantran2012,valizadehnatarajan2013} and a very good agreement can be observed. Next, we illustrate the performance of the present formulation for thin plate problems. A simply supported square plate subjected to uniform load is considered, while the length-to-thickness $(a/h)$ varies from 5 to 10$^4$. Two individual approaches are employed: discrete shear gap method referred to as DSG3 and the other referred as cell-based smoothed finite element method with discrete shear gap technique (CSDSG3). The normalized center deflection $\overline{w}_c = 100 w_c \frac{E_ch^3}{12(1-\nu^2)pa^4}$ is shown in \fref{fig:plateshearlock}. It is observed that the DSG3 results are subjected to shear locking when the plate becomes thin $(a/h > 100)$. However, the present formulation, CSDSG3 is less sensitive to shear locking. 

\begin{table}[htbp]
\centering
\renewcommand{\arraystretch}{1.2}
\caption{The normalized center deflection $\overline{w}_c = 100 w_c \frac{E_ch^3}{12(1-\nu^2)pa^4}$ for a simply supported Al/ZrO$_2$-1 FGM square plate with $a/h=$ 5, subjected to a uniformly distributed load $p$.}
\begin{tabular}{llccc}
\hline
Method & \multicolumn{3}{c}{gradient index, $n$} \\
\cline{2-4}
&  0 & 1 & 2 \\
\hline
4$\times$4&0.1443	&0.2356&	0.2644\\
8$\times$8	&0.1648	&0.2703&0.3029\\
16$\times$16&	0.1701	&0.2795&	0.3131\\
32$\times$32	&0.1714&	0.2819&	0.3158\\
40$\times$40	&0.1716&	0.2822	& 0.3161 \\
NS-DSG3~\cite{nguyen-xuantran2012} & 0.1721  & 0.2716 & 0.3107 \\
ES-DSG3~\cite{nguyen-xuantran2012} & 0.1700 & 0.2680 & 0.3066 \\
MLPG~\cite{gilhooleybatra2007} & 0.1671  & 0.2905 & 0.3280 \\
$kp-$Ritz~\cite{leezhao2009} & 0.1722  & 0.2811 & 0.3221 \\
MITC4~\cite{nguyen-xuantran2012} & 0.1715  & 0.2704 & 0.3093 \\
IGA-Quadratic~\cite{valizadehnatarajan2013} & 0.1717  & 0.2719 & 0.3115 \\
\hline
\end{tabular}
\label{tab:ssss-valid}
\end{table}

\subsection{Free flexural vibrations}
\begin{figure}[htpb]
\centering
\includegraphics[scale=0.7]{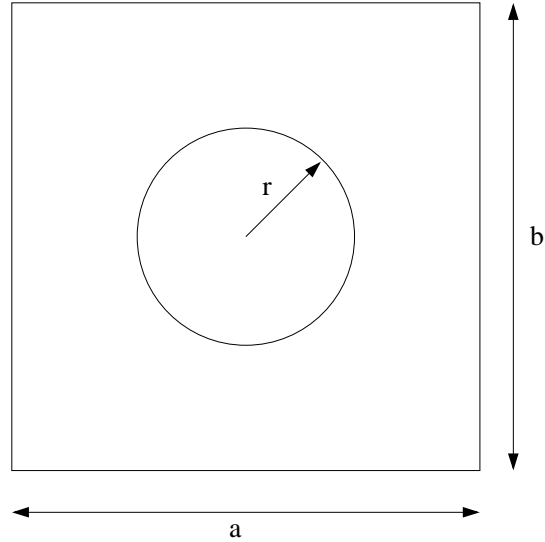}
\caption{Plate with a centrally located circular cutout. $r$ is the radius of the circular cutout.}
\label{fig:phole}
\end{figure}

In this section, the free flexural vibration characteristics of FGM plates with and without centrally located cutout in thermal environment is studied numerically. \fref{fig:phole} shows the geometry of the plate with a centrally located circular cutout. In all cases, we present the non-dimensionalized free flexural frequency defined as, unless otherwise stated:

\begin{equation}
\overline{\omega} = \omega a^2 \sqrt{ \frac{\rho_c h}{D_c}}
\end{equation}

\begin{figure}[htpb]
\centering
\includegraphics[scale=0.7]{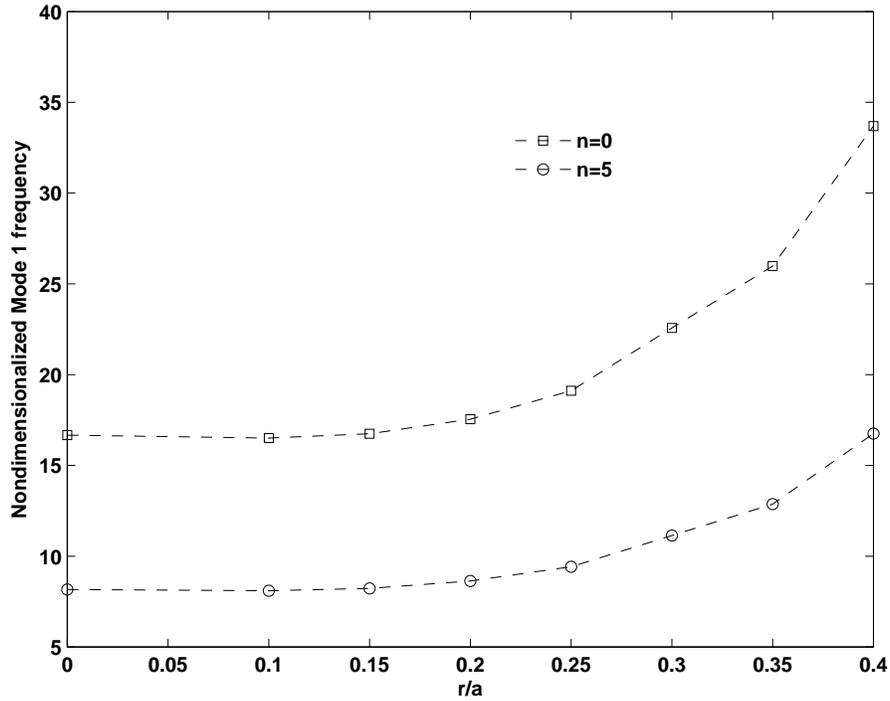}
\caption{Effect of the cutout size on the fundamental frequency $(\Omega)$ for a square simply supported FGM plate with a central circular cutout in thermal environment $\Delta T =$ 100K $(T_c =$ 400K, $T_m=$ 300K$)$ for different gradient index $n$.}
\label{fig:vibcutoutsize}
\end{figure}

where $\omega$ is the natural frequency, $\rho_c, D_c = \frac{E_c h^3}{12(1-\nu^2)}$ are the mass density and the flexural rigidity of the ceramic phase. The FGM plate considered here is made up of silicon nitride (Si$_3$N$_4$) and stainless steel (SUS304). The material is considered to be temperature dependent and the temperature coefficients corresponding to Si$_3$N$_4$/SUS304 are listed in Table \ref{table:tempdepprop}~\cite{Reddy1998,Sundararajan2005}. The mass density $(\rho)$ and the thermal conductivity $(\kappa)$ are $\rho_c=$ 2370 kg/m$^3$, $\kappa_c=$ 9.19 W/mK for Si$_3$N$_4$ and $\rho_m=$ 8166 kg/m$^3$, $\kappa_m=$ 12.04 W/mK for SUS304. Poisson's ratio $\nu$ is assumed to be constant and taken as 0.28 for the current study~\cite{Sundararajan2005}. Before proceeding with a detailed study on the effect of gradient index on the natural frequencies, the formulation developed herein is validated against available analytical/numerical solutions pertaining to the linear frequencies of a FGM plate in thermal environment and a FGM plate with a centrally located circular cutout. The computed frequencies: (a) for a square simply supported FGM plate in thermal environment with $a/h=$ 10 is given in Table \ref{tab:SSfgmValida} and (b) the mesh convergence and comparison of linear frequencies for a square plate with circular cutout is given in Tables \ref{table:vibcutconvergence} - \ref{table:vibgIndexcutoutvalida}. It can be seen that the numerical results from the present formulation are found to be in very good agreement with the existing solutions. For the uniform temperature case, the material properties are evaluated at $T_c = T_m=$ 300K. The temperature is assumed vary only in the thickness direction and determined by \Eref{eqn:tempsolu}. The temperature for the ceramic surface is varied, whist maintaining a constant value on the metallic surface is maintained $(T_m=$ 300K$)$ to subject a thermal gradient. The geometric stiffness matrix is computed from the in-plane stress resultants due to the applied thermal gradient. The geometric stiffness matrix is then added to the stiffness matrix and the eigenvalue problem is solved. The effect of the material gradient index is also shown in Tables \ref{tab:SSfgmValida} \& \ref{table:vibgIndexcutoutvalida} and the influence of a centrally located cutout is shown in Tables \ref{table:vibcutconvergence} - \ref{table:vibgIndexcutoutvalida}. The combined effect of increasing the temperature and the gradient index is to lower the fundamental frequency, this is due to the increase in the metallic volume fraction. \fref{fig:vibcutoutsize} shows the influence of the cutout size on the frequency for a plate in thermal environment $(\Delta T =$ 100K$)$. The frequency increases with increasing cutout size. This can be attributed to the decrease in stiffness due to the presence of the cutout. Also, it can be seen that with increasing gradient index, the frequency decreases. In this case, the decrease in the frequency is due to the increase in the metallic volume fraction. It is observed that the combined effect of increasing the gradient index and the cutout size is to lower the fundamental frequency. Increasing the thermal gradient further decreases the fundamental frequency.

\begin{table}
\renewcommand\arraystretch{1.5}
\caption{Temperature dependent coefficient for material Si$_3$N$_4$/SUS304, Ref~\cite{Reddy1998,Sundararajan2005}.}
\centering
\begin{tabular}{lcccccc}
\hline
Material & Property & $P_o$ & $P_{-1}$ & $P_1$ & $P_2$ & $P_3$  \\
\hline
\multirow{2}{*}{Si$_3$N$_4$} & $E$(Pa) & 348.43e$^9$ &0.0& -3.070e$^{-4}$ & 2.160e$^{-7}$ & -8.946$e^{-11}$  \\
& $\alpha$ (1/K) & 5.8723e$^{-6}$ & 0.0 & 9.095e$^{-4}$ & 0.0 & 0.0 \\
\cline{2-7}
\multirow{2}{*}{SUS304} & $E$(Pa) & 201.04e$^9$ &0.0& 3.079e$^{-4}$ & -6.534e$^{-7}$ & 0.0  \\
& $\alpha$ (1/K) & 12.330e$^{-6}$ & 0.0 & 8.086e$^{-4}$ & 0.0 & 0.0 \\
\hline
\end{tabular}
\label{table:tempdepprop}
\end{table}

\begin{table}[htbp]
\centering
\renewcommand{\arraystretch}{1.5}
\caption{The first normalized frequency parameter $\overline{\omega} $ for a fully simply supported Si$_3$N$_4$/SUS304 FGM square plate with $a/h=$ 10 in thermal environment.}
\begin{tabular}{llcccc}
\hline
$T_c$,$T_m$ & &\multicolumn{4}{c}{gradient index $n$} \\
\cline{3-6}
 & & 0 & 1 & 5 & 10 \\
\hline
\multirow{2}{*}{300,300}& Present & 18.3570 & 11.0690 & 9.0260 & 8.5880 \\
& Ref.~\cite{natarajanbaiz2011b} & 18.3731 & 11.0288 & 9.0128 & 8.5870 \\
\multirow{2}{*}{400,300}& Present & 17.9778 & 10.7979 & 8.8626 & 8.3182 \\
& Ref.~\cite{natarajanbaiz2011b} & 17.9620 & 10.7860 & 8.7530 & 8.3090 \\
\multirow{2}{*}{600,300} & Present & 17.1205 & 10.1679 & 8.1253 & 7.6516 \\
& Ref.~\cite{natarajanbaiz2011b} & 17.1050 & 10.1550 & 8.1150 & 7.6420 \\ 
\hline
\end{tabular}%
\label{tab:SSfgmValida}%
\end{table}%

\begin{table}[htpb]
\renewcommand\arraystretch{1.}
\caption{Convergence of fundamental frequency $ \left( \Omega = \left[ \frac{\omega^2 \rho_c h a^4}{D_c(1-\nu^2)} \right]^{1/4} \right) $ with mesh size for an isotropic plate with a central cutout.}
\centering
\begin{tabular}{lcc}
\hline 
Number of nodes & Mode 1 & Mode 1 \\
\hline
333 & 6.1025 & 8.6297 \\
480 & 6.0805 & 8.5595 \\
719 & 6.0663 & 8.5192 \\
1271 & 6.0560	 & 8.4852 \\
Ref.~\cite{ahmadnatarajan2013} & 6.1725 & 8.6443 \\
Ref.~\cite{huangsakiyama1999} & 6.2110 & 8.7310 \\
\hline
\end{tabular}
\label{table:vibcutconvergence}
\end{table}

\begin{table}[htpb]
\renewcommand\arraystretch{1.5}
\caption{Comparison of fundamental frequency for a simply supported FGM plate with $a/h=$ 5 and $r/a =$ 0.2. }
\centering
\begin{tabular}{clccccc}
\hline 
$T_c$ & & \multicolumn{5}{c}{gradient index, $n$} \\
\cline{3-7}
 &  & 0 & 1 & 2 & 5 & 10 \\
 \hline
 \multirow{2}{*}{300} & Ref.~\cite{ahmadnatarajan2013} & 17.6855 & 10.6681 & 9.6040 & 8.7113 & 8.2850 \\
 & Present & 17.7122 & 10.6845 & 9.6188 & 8.7246 & 8.2976 \\
 \multirow{2}{*}{400} & Ref.~\cite{ahmadnatarajan2013} & 17.4690 & 10.5174 & 9.4618 & 8.5738 & 8.1484 \\
 & Present & 17.5488 & 10.5775 & 9.5197 & 8.6309 & 8.2059 \\
\hline
\end{tabular}
\label{table:vibgIndexcutoutvalida}
\end{table}

\subsection{Buckling analysis}
In this section, we present the mechanical and thermal buckling behaviour of functionally graded skew plates. 

\subsubsection*{Mechanical Buckling}
\begin{figure}[htpb]
\centering
\includegraphics[scale=0.7]{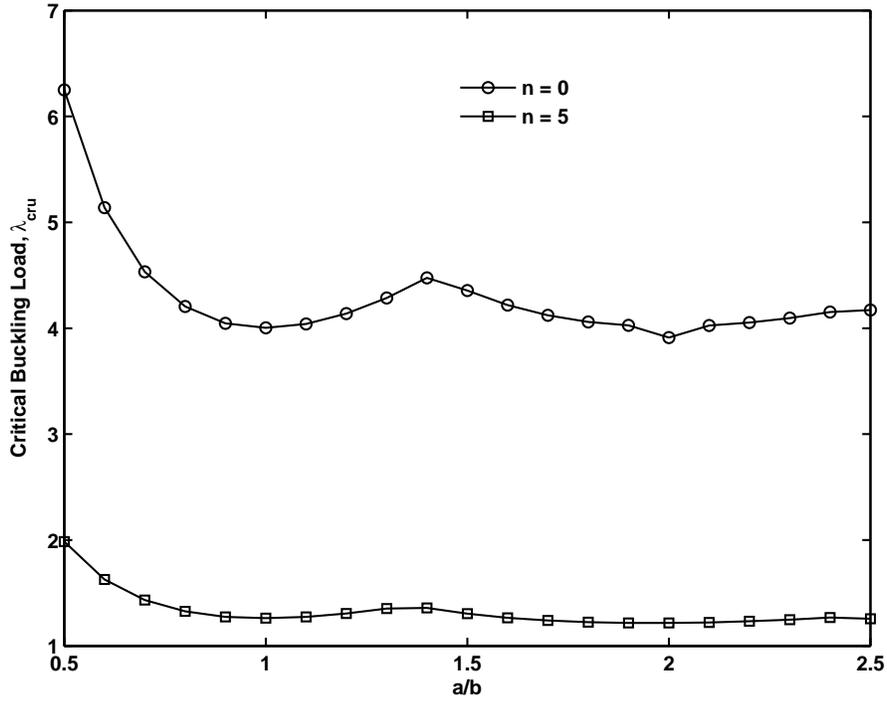}
\caption{Effect of plate aspect ratio $a/b$ and gradient index on the critical buckling load for a simply supported FGM plate under uni-axial compression with $a/h=10$.}
\label{fig:mbuck_ab}
\end{figure}

\begin{figure}[htpb]
\centering
\includegraphics[scale=0.7]{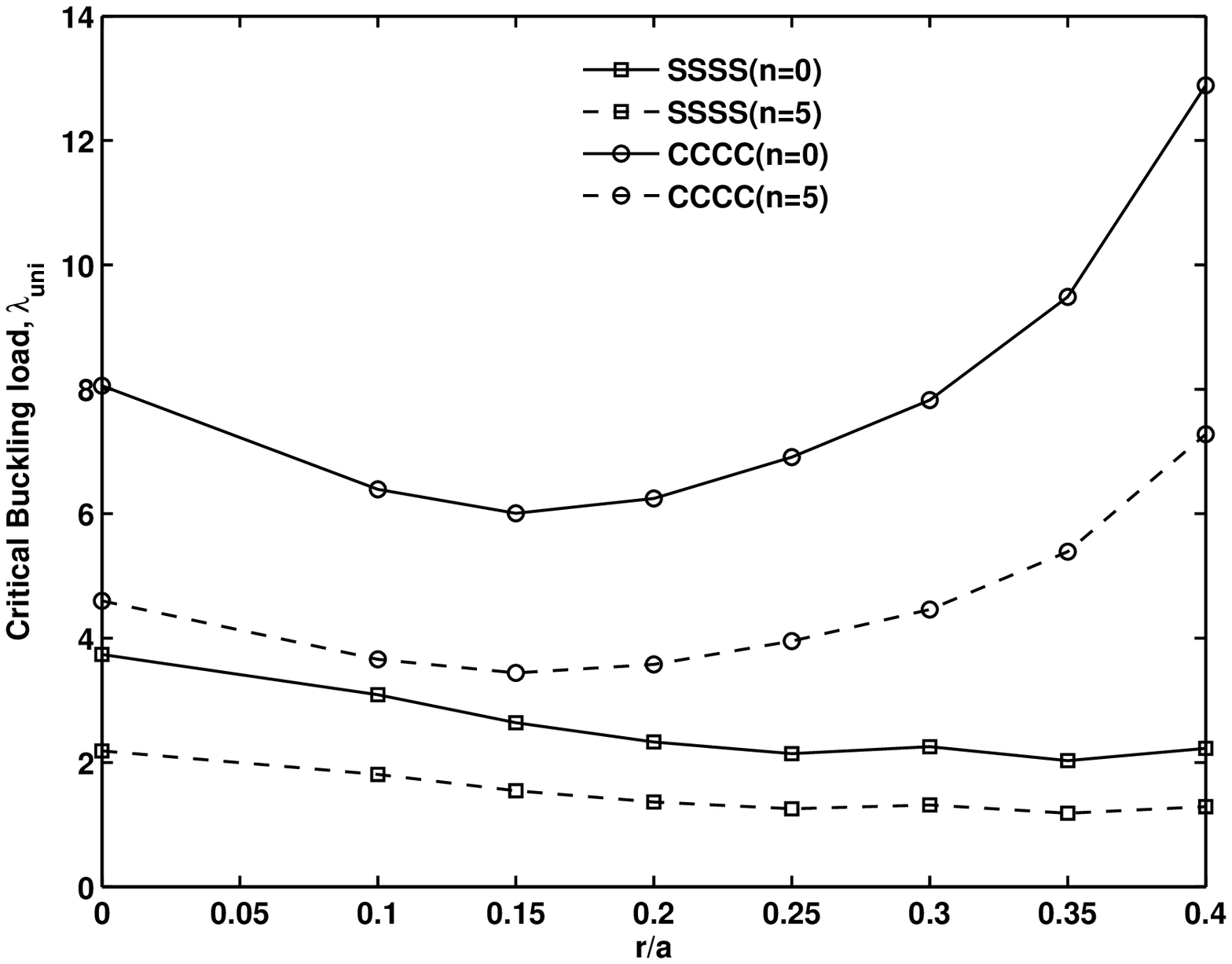}
\caption{Variation of the critical buckling load, $\lambda_{cru} = \frac{N_{xxcr}^o b^2}{\pi^2 D_c}$ with cutout dimensions for a square FGM plate with a/h=10 subjected to uniaxial compressive loading for different gradient index $n$ and various boundary conditions.}
\label{fig:mbuckravariation}
\end{figure}

The FGM plate considered here consists of Aluminum (Al) and Zirconium dioxide (ZrO$_2$). The material is considered to be temperature independent. The Young's modulus $(E)$ for ZrO$_2$ is $E_c=$ 151 GPa and for Al is $E_m=$ 70 GPa. For mechanical buckling, we consider both uni- and bi-axial mechanical loads on the FGM plates. In all cases, we present the critical buckling parameters as, unless otherwise specified:

\begin{eqnarray}
\lambda_{\rm{cru}} = \frac{N_{\rm{xxcr}}^0 b^2}{\pi^2 D_c} \nonumber \\
\lambda_{\rm{crb}} = \frac{N_{\rm{yycr}}^0 b^2}{\pi^2 D_c} 
\label{eqn:mbuckparm}
\end{eqnarray}
where $\lambda_{\rm{cru}}$ and $\lambda_{\rm{crb}}$ are the critical buckling parameters for uni- and bi-axial load, respectively, $D_c=E_ch^3/(12(1-\nu^2))$. The critical buckling loads evaluated by varying the skew angle of the plate, volume fraction index and considering mechanical loads (uni- and biaxial compressive loads) are shown in Tables \ref{table:mbuckskewangleah100} for $a/h=$ 100. The efficacy of the present formulation is demonstrated by comparing our results with those in~\cite{ganapathiprakash2006}. It can be seen that increasing the gradient index decreases the critical buckling load. A very good agreement in the results can be observed. It is also observed that the decrease in the critical value is significant for the material gradient index $n \le 2$ and that further increase in $n$ yields less reduction in the critical value, irrespective of the skew angle. The effect of the plate aspect ratio and the gradient index on the critical buckling load is shown in \fref{fig:mbuck_ab} for a simply supported FGM plate under uni-axial mechanical load. It is observed that the combined effect of increasing the gradient index and the plate aspect ratio is to lower the critical buckling load. Table \ref{table:mbuckgIndex} presents the critical buckling parameter for a simply supported FGM with a centrally located circular cutout with $r/a=$ 0.2. It can be seen that the present formulation yields comparable results. The effect of increasing the gradient index is to lower the critical buckling load. This can be attributed to the stiffness degradation due to increase in the metallic volume fraction. \fref{fig:mbuckravariation} shows the influence of a centrally located circular cutout and the gradient index on the critical buckling load under two different boundary conditions, viz., all edges simply supported and all edges clamped. In this case, the plate is subjected to a uni-axial compressive load. It can be seen that increasing the gradient index decreases the critical buckling load due to increasing metallic volume fraction, whilst, increasing the cutout radius decreases the critical buckling load in the case of simply supported boundary conditions. This can be attributed to the stiffness degradation due to the presence of a cutout and the simply supported boundary condition. In case of the clamped boundary condition, the critical buckling load first decreases with increasing cutout radius due to stiffness degradation. Upon further increase, the critical buckling load increases. This is because, the clamped boundary condition adds stiffness to the system which overcomes the stiffness reduction due to the presence of a cutout. 

\begin{table}[htpb]
\renewcommand\arraystretch{1.5}
\caption{Critical buckling parameters for a thin simply supported FGM skew plate with $a/h=$ 100 and $a/b=$ 1.}
\centering
\begin{tabular}{ccccccccc}
\hline 
Skew angle & $\lambda_{cr}$ & \multicolumn{6}{c}{Gradient index, $n$} \\
\cline{3-8}
 &  & \multicolumn{2}{c}{0} & \multicolumn{2}{c}{1} &  5 & 10\\
\cline{3-6}
 &  & Ref.~\cite{ganapathiprakash2006} & Present & Ref.~\cite{ganapathiprakash2006} & Present \\
\hline
\multirow{2}{*}{0$^\circ$} & $\lambda_{\rm{cru}}$ & 4.0010 & 4.0034 & 1.7956 & 1.8052 & 1.2624 & 1.0846 \\
& $\lambda_{\rm{crb}}$ & 2.0002 & 2.0017 & 0.8980 & 0.9028 &  0.6312 & 0.5423 \\
\multirow{2}{*}{15$^\circ$} & $\lambda_{\rm{cru}}$ & 4.3946 & 4.4007 & 1.9716 & 1.9799 & 1.3859 & 1.1915 \\
& $\lambda_{\rm{crb}}$ & 2.1154 & 2.1187 & 0.9517 & 0.9561 &  0.6683 & 0.5741 \\
\multirow{2}{*}{30$^\circ$} & $\lambda_{\rm{cru}}$ & 5.8966 & 5.9317 & 2.6496 & 2.6496 & 1.8586 & 1.6020 \\
& $\lambda_{\rm{crb}}$ & 2.5365 & 2.5491 & 1.1519 & 1.1520 & 0.8047 & 0.6909 \\
\hline
\end{tabular}
\label{table:mbuckskewangleah100}
\end{table}

\begin{table}[htpb]
\renewcommand\arraystretch{1.5}
\caption{Comparison of critical buckling load $\lambda_{cru} = \frac{N_{xxcr}^o b^2}{\pi^2 D_m}$ for a simply supported FGM plate with $a/h=$ 100 and $r/a =$ 0.2. The effective material properties are computed by rule of mixtures. In order to be consistent with the literature, the properties of the metallic phase is used for normalization.}
\centering
\begin{tabular}{lccc}
\hline 
gradient index, $n$ & Ref.~\cite{zhaolee2009a} & Present & \% difference \\
\hline
0 & 5.2611 & 5.2831 & -0.42 \\
0.2 & 4.6564 & 4.6919 & -0.76 \\
1 & 3.6761 & 3.663 & 0.36 \\
2 & 3.3672 & 3.3961 & -0.86 \\
5 & 3.1238 & 3.1073 & 0.53 \\
10 & 2.9366	 & 2.8947 & 1.43 \\
\hline
\end{tabular}
\label{table:mbuckgIndex}
\end{table}

\subsubsection*{Thermal Buckling}
\begin{figure}[htpb]
\centering
\includegraphics[scale=0.7]{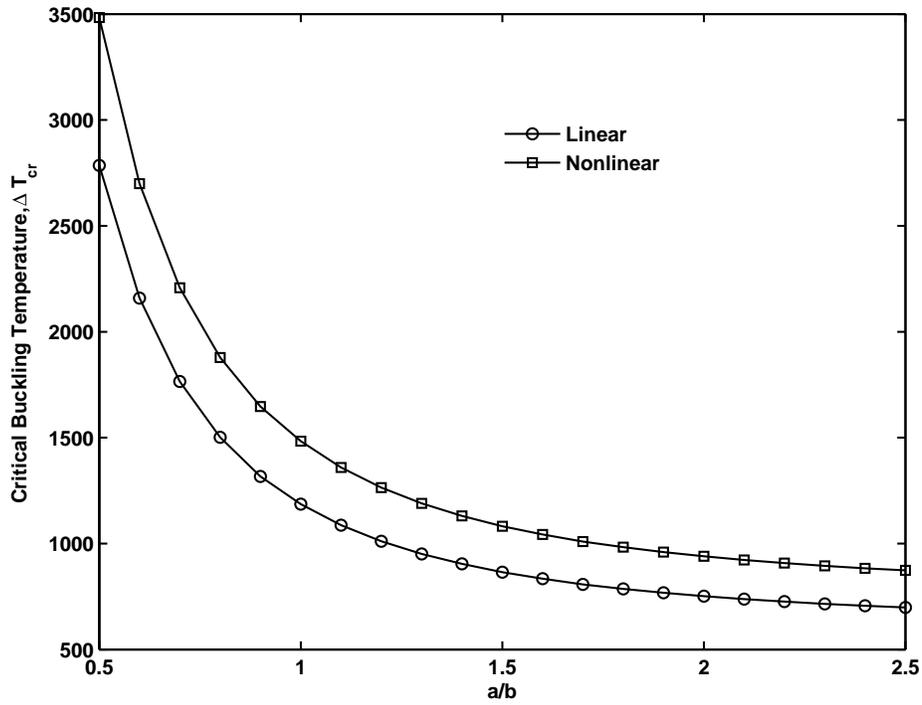}
\caption{Critical buckling temperature as a function of plate aspect ratio $a/b$ with linear and nonlinear temperature distribution through the thickness with gradient index $n=$ 5.}
\label{fig:tbuck_ab}
\end{figure}

\begin{figure}[htpb]
\centering
\includegraphics[scale=0.7]{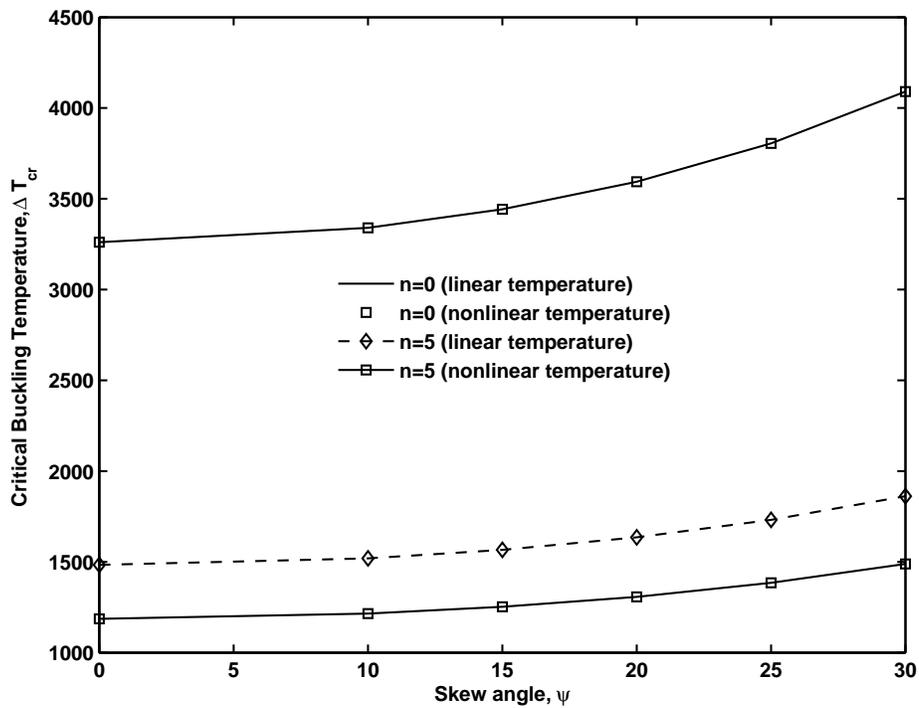}
\caption{Critical buckling temperature as a function of skew angle $\psi$ for a simply supported square FGM plate with a/h=10. Both linear and nonlinear temperature distribution through the thickness is assumed.}
\label{fig:tbuck_skew}
\end{figure}

\begin{figure}[htpb]
\centering
\includegraphics[scale=0.7]{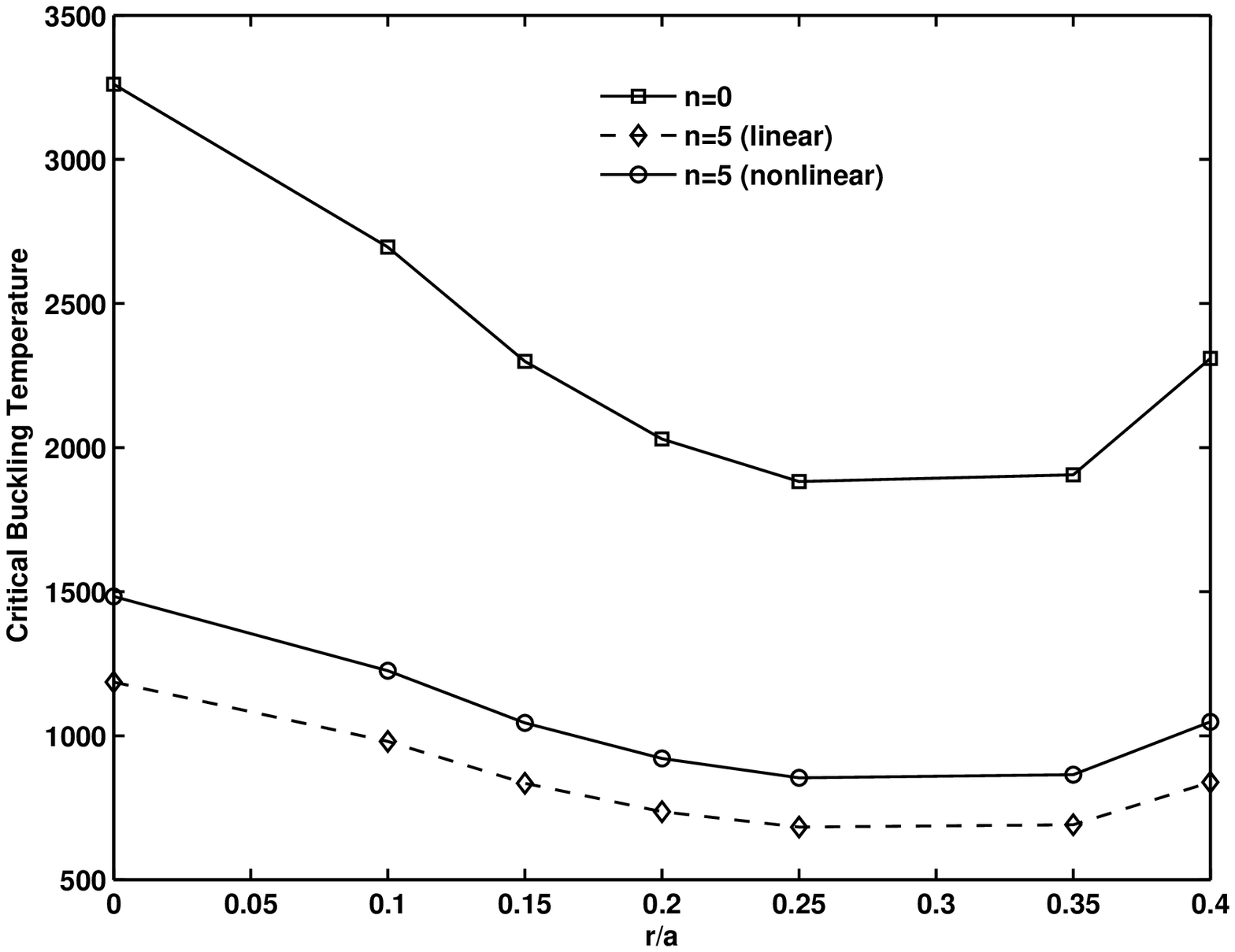}
\caption{Influence of cutout size on the critical buckling temperature for a square simply supported FGM plate with a centrally located circular cutout with $a/h=$ 10 for various gradient index $n$. Linear and nonlinear temperature distribution through the thickness is assumed.}
\label{fig:ThBuckRa}
\end{figure}

The thermal buckling behaviour of simply supported functionally graded skew plate is studied next. The top surface is ceramic rich and the bottom surface is metal rich. The FGM plate considered here consists of aluminum and alumina. The Young's modulus, the thermal conductivity and the coefficient of thermal expansion for alumina is $E_c=$ 380 GPa, $K_c=$10.4 W/mK, $\alpha_c=$ 7.4 $\times$ 10$^{-6}$ 1/$^\circ$C, and for aluminum, $E_m=$ 70 GPa, $K_m=$ 204 W/mK, $\alpha_m=$ 23 $\times$ 10$^{-6}$ 1/$^\circ$C, respectively. Poisson's ratio is chosen as constant, $\nu=$ 0.3. The temperature rise of $T_m=$ 5$^\circ$C in the metal-rich surface of the plate is assumed in the present study. In addition to nonlinear temperature distribution across the plate thickness, the linear case is also considered in the present analysis by truncating the higher order terms in \Eref{eqn:heatconducres}. The plate is of uniform thickness and simply supported on all four edges. Table \ref{table:tbuckvalid} shows the convergence of the critical buckling temperature with mesh size for different gradient index, $n$. It can be seen that the results from the present formulation are in good very agreement with the available solution. The influence of the plate aspect ratio $a/b$ and the skew angle $\psi$ on the critical buckling temperature for a simply supported square FGM plates are shown in \frefs{fig:tbuck_ab} and (\ref{fig:tbuck_skew}). It is seen that increasing the plate aspect ratio decreases the critical buckling temperature for both linear and nonlinear temperature distribution through the thickness. The critical buckling temperature increases with increase in the skew angle. The influence of the gradient index $n$ is also shown in \fref{fig:tbuck_skew}. It is seen that with increasing gradient index, $n$, the critical buckling temperature decreases. This is due to the increase in the metallic volume fraction that degrades the overall stiffness of the structure. \fref{fig:ThBuckRa} shows the influence of the cutout radius and the material gradient index on the critical buckling temperature. Both linear and nonlinear temperature distribution through the thickness is assumed. Again, it is seen that the combined effect of increasing the gradient index $n$ and the cutout radius $r/a$ is to lower the buckling temperature. For gradient index $n=$ 0, there is no difference between the linear and the nonlinear temperature distribution through the thickness as the material is homogeneous through the thickness. While, for $n > 0$, the material is heterogeneous through the thickness with different thermal property.

\begin{table}[htpb]
\renewcommand\arraystretch{1.5}
\caption{Convergence of the critical buckling temperature for a simply supported FGM skew plate with $a/h=$ 10 and $a/b=$ 1. Nonlinear temperature rise through the thickness of the plate is assumed.}
\centering
\begin{tabular}{lcccc}
\hline 
Mesh & \multicolumn{4}{c}{Gradient index, $n$} \\
\cline{2-5}
 & 0 & 1 & 5 & 10 \\
\hline
8$\times$8 & 3383.40 & 2054.61 & 1539.24 & 1496.36 \\
16$\times$16 & 3286.90 & 1995.07 & 1495.25 & 1453.99 \\
32$\times$32 & 3263.91 & 1980.96 & 1484.76 & 1443.86 \\
40$\times$40 & 3261.17 & 1979.30 & 1483.51 & 1442.60 \\		
Ref.~\cite{ganapathiprakash2006a} & 3257.47 & 1977.01 & 1481.83 & 1441.02 \\
\hline
\end{tabular}
\label{table:tbuckvalid}
\end{table}

\section{Conclusion}
In this paper, we applied the cell-based smoothed finite element  method with discrete shear gap technique to study the static and the dynamic response of functionally graded materials. The first order shear deformation theory was used to describe the plate kinematics. The efficiency and accuracy of the present approach is demonstrated with few numerical examples. This improved finite element technique shows insensitivity to shear locking and produce excellent results in static bending, free vibration and buckling of functionally graded plates. 

\section*{Acknowledgements} 
S Natarajan would like to acknowledge the financial support of the School of Civil and Environmental Engineering, The University of New South Wales for his research fellowship since September 2012.

\newpage
\section*{References}
\bibliographystyle{plain}
\bibliography{laminatedCompo}

\begin{thebibliography}{10}

\bibitem{bletzingerbischoff2000}
KU~Bletzinger, M~Bischoff, and E~Ramm.
\newblock A unified approach for shear locking free triangular and rectangular
  shell finite elements.
\newblock {\em Computers and Structures}, 75:321--334, 2000.

\bibitem{bordasnatarajan2010}
S~Bordas and S~Natarajan.
\newblock On the approximation in the smoothed finite element method {(SFEM)}.
\newblock {\em International Journal for Numerical Methods in Engineering},
  81:660--670, 2010.

\bibitem{bordasnatarajan2011}
S~Bordas, S~Natarajan, P~Kerfriden, CE~Augarde, D~Roy Mahapatra, T~Rabczuk, and
  S~Dal Pont.
\newblock {On the performance of strain smoothing for quadratic and enriched
  finite el- ement approximation (XFEM/GFEM/PUFEM)}.
\newblock {\em International Journal for Numerical Methods in Engineering},
  86:637--666, 2011.

\bibitem{carrrerabrischetto2011}
E~Carrera, S~Brischetto, and M~Cinefra.
\newblock Effects of thickness stretching in functionally graded plates and
  shells.
\newblock {\em Composites Part N: Engineering}, 42:123--133, 2011.

\bibitem{carrerabrischetto2008}
E~Carrera, S~Brischetto, and A~Robaldo.
\newblock Variable kinematic model for the analysis of functionally graded
  material plates.
\newblock {\em AIAA Journal}, 46:194--203, 2008.

\bibitem{cazesmeschke2012}
F~Cazes and G~Meschke.
\newblock An edge based imbricate finite element method {(EI-FEM)} with full
  and reduced integration.
\newblock {\em Computer \& Structures}, 106--107:154--175, 2012.

\bibitem{cinefracarrera2012}
M~Cinefra, E~Carrera, L~Della Croce, and C~Chinosi.
\newblock Refined shell elements for the analysis of functionally graded
  structures.
\newblock {\em Composite Structures}, 94:415--422, 2012.

\bibitem{dailiu2005}
KY~Dai, GR~Liu, X~Han, and KM~Lim.
\newblock Thermomechanical analysis of functionally graded material (fgm)
  plates using element-free galerkin method.
\newblock {\em Computers and Structures.}, 83:1487--1502, 2011.

\bibitem{ferreirabatra2006}
AJM Ferreira, RC~Batra, CMC Roque, LF~Qian, and RMN Jorge.
\newblock Natural frequencies of functionally graded plates by a meshless
  method.
\newblock {\em Composite Structures}, 75:593--600, 2006.

\bibitem{ganapathiprakash2006a}
M~Ganapathi and T~Prakash.
\newblock Thermal buckling of simply supported functionally graded skew plates.
\newblock {\em Composite Structures}, 74:247--250, 2006.

\bibitem{ganapathiprakash2006}
M~Ganapathi, T~Prakash, and N~Sundararajan.
\newblock Influence of functionally graded material on buckling of skew plates
  under mechanical loads.
\newblock {\em Journal of Engineering Mechanics: ASCE}, 1332:902--905, 2006.

\bibitem{gilhooleybatra2007}
DF~Gilhooley, RC~Batra, JR~Xiao, MA~McCarthy, and JW~Gillespie.
\newblock Analysis of thick functionally graded plates by using higher order
  shear and normal deformable plate theory and {MLPG} method with radial basis
  functions.
\newblock {\em Composite Structures}, 80:539--552, 2007.

\bibitem{hecheng2011}
ZC~He, AG~Cheng, GY~Zhang, ZH~Zhong, and GR~Liu.
\newblock Dispersion error reduction for acoustic problems using the edge based
  smoothed finite element method {(ES-FEM)}.
\newblock {\em International Journal for Numerical Methods in Engineering},
  86:1322--1338, 2011.

\bibitem{huangsakiyama1999}
M~Huang and T~Sakiyama.
\newblock Free vibration analysis of rectangular plates with variously shaped
  holes.
\newblock {\em Journal of Sound and Vibration}, 226:769--786, 1999.

\bibitem{janghorbanazare2011}
M~Janghorbana and A~Zare.
\newblock Thermal effect on free vibration analysis of functionally graded
  arbitrary straight-sided plates with different cutouts.
\newblock {\em Latin American Journal of Solids and Structures}, 8:245--257,
  2011.

\bibitem{jhakant2013}
DK~Jha, Tarun Kant, and RK~Singh.
\newblock A critical review of recent research on functionally graded plates.
\newblock {\em Composite Structures}, 96:833--849, 2013.

\bibitem{koizumi1993}
M~Koizumi.
\newblock The concept of {FGM}.
\newblock {\em Ceramic Transactions - Functionally graded materials}, 34:3--10,
  1993.

\bibitem{leezhao2009}
YY~Lee, X~Zhao, and KM~Liew.
\newblock Thermo-elastic analysis of functionally graded plates using the
  element free $kp-${R}itz method.
\newblock {\em Smart Materials and Structures}, 18:035007, 2009.

\bibitem{liunguyen2008}
G~Liu, T~Nguyen-Thoi, and K~Lam.
\newblock A novel alpha finite element method ($\alpha$fem) for exact solution
  to mechanics problems using triangular and tetrahedral elements.
\newblock {\em Computer Methods in Applied Mechanics and Engineering},
  197:3883--3897, 2008.

\bibitem{liunguyen2009a}
G~Liu, T~Nguyen-Thoi, and K~Lam.
\newblock An edge-based smoothed finite element method ({ES-FEM}) for static,
  free and forced vibration analyses of solids.
\newblock {\em Journal of Sound and Vibration}, 320:1100--1130, 2009.

\bibitem{liunguyen2009}
G~Liu, T~Nguyen-Thoi, H~Nguyen-Xuan, and K~Lam.
\newblock A node based smoothed finite element ({NS-FEM}) for upper bound
  solution to solid mechanics problems.
\newblock {\em Computers and Structures}, 87:14--26, 2009.

\bibitem{liudai2007}
GR~Liu, KY~Dai, and TT~Nguyen.
\newblock A smoothed finite element for mechanics problems.
\newblock {\em Computational Mechanics}, 39:859--877, 2007.

\bibitem{natarajanbaiz2011}
S~Natarajan, Pedro~M Baiz, S~Bordas, T~Rabczuk, and P~Kerfriden.
\newblock Natural frequencies of cracked functionally graded material plates by
  the extended finite element method.
\newblock {\em Composite Structures}, 93:3082--3092, 2011.

\bibitem{natarajanbaiz2011b}
S~Natarajan, PM~Baiz, M~Ganapathi, P~Kerfriden, and S~Bordas.
\newblock Linear free flexural vibration of cracked functionally graded plates
  in thermal environment.
\newblock {\em Computers and Structures}, 89:1535--1546, 2011.

\bibitem{thoiliu2009}
T~Nguyen-Thoi, G~Liu, K~Lam, and G~Zhang.
\newblock A face-based smoothed finite element method ({FS-FEM}) for 3{D}
  linear and nonlinear solid mechanics using 4-node tetrahedral elements.
\newblock {\em International Journal for Numerical Methods in Engineering},
  78:324--353, 2009.

\bibitem{nguyenbordas2008}
H~Nguyen-Xuan, S~Bordas, and H~Nguyen-Dang.
\newblock Smooth finite element methods: convergence, accuracy and properties.
\newblock {\em International Journal for Numerical Methods in Engineering},
  74:175--208, 2008.

\bibitem{nguyenliu2013}
H~Nguyen-Xuan, GR~Liu, S~Bordas, S~Natarajan, and T~Rabczuk.
\newblock Ad adaptive singular {ES-FEM} for mechanics problems with singular
  field of arbitrary order.
\newblock {\em Computer Methods in Applied Mechanics and Engineering},
  253:252--273, 2013.

\bibitem{nguyenrabczuk2008}
H~Nguyen-Xuan, T~Rabczuk, S~Bordas, and JF~Debongnie.
\newblock A smoothed finite element method for plate analysis.
\newblock {\em Computer Methods in Applied Mechanics and Engineering},
  197:1184--1203, 2008.

\bibitem{nguyen-xuantran2012}
H~Nguyen-Xuan, Loc~V Tran, H~Thai, and T~Nguyen-Thoi.
\newblock Analysis of functionally graded plates by an efficient finite element
  method with node-based strain smoothing.
\newblock {\em Thin Walled Structures}, 54:1--18, 2012.

\bibitem{nguyenrabczuk2008a}
NT~Ngyuen, T~Rabczuk, H~Nguyen-Xuan, and S~Bordas.
\newblock A smoothed finite element method for shell analysis.
\newblock {\em Computer Methods in Applied Mechanics and Engineering},
  198:165--177, 2008.

\bibitem{praveenreddy1998}
GN~Praveen and JN~Reddy.
\newblock Nonlinear transient thermoelastic ceramic-metal plates.
\newblock {\em International Journal of Solids and Structures}, 35:4457--4476,
  1998.

\bibitem{qianbatra2004}
LC~Qian, RC~Batra, and LM~Chen.
\newblock Static and dynamic deformations of thick functionally graded elastic
  plates by using higher order shear and normal deformable plate theory and
  meshless local {Petrov Galerkin} method.
\newblock {\em Composites Part B: Engineering}, 35:685--697, 2004.

\bibitem{ahmadnatarajan2013}
AA~Rahimabadi, S~Natarajan, and S~Bordas.
\newblock Vibration of functionally graded material plates with cutouts \&
  cracks in thermal environment.
\newblock {\em Key Engineering Materials}, 560:157--180, 2013.

\bibitem{Rajasekaran1973}
S~Rajasekaran and DW~Murray.
\newblock Incremental finite element matrices.
\newblock {\em {ASCE} Journal of Structural Divison}, 99:2423--2438, 1973.

\bibitem{reddy2000}
J.~N. Reddy.
\newblock Analysis of functionally graded plates.
\newblock {\em International Journal for Numerical Methods in Engineering},
  47:663--684, 2000.

\bibitem{Reddy1998}
JN~Reddy and CD~Chin.
\newblock Thermomechanical analysis of functionally graded cylinders and
  plates.
\newblock {\em Journal of Thermal Stresses}, 21:593--629, 1998.

\bibitem{reddychin2007}
JN~Reddy and CD~Chin.
\newblock Thermomechanical analysis of functionally graded cylinders and
  plates.
\newblock {\em Journal of Thermal Stresses}, 21:593--626, 2007.

\bibitem{sajivarughese2008}
D~Saji, Byji Varughese, and SC~Pradhan.
\newblock Finite element analysis for thermal buckling behaviour in
  functionally graded plates with cutouts.
\newblock volume 113, 2008.

\bibitem{shariateslami2006}
BA~Samsam Shariat and MR~Eslami.
\newblock Thermal buckling of imperfect functionally graded plates.
\newblock {\em International Journal of Solids and Structures.}, 43:4082--4096,
  2006.

\bibitem{singhaprakash2011}
MK~Singha, T~Prakash, and M~Ganapathi.
\newblock Finite element analysis of functionally graded plates.
\newblock {\em Finite Elements in Analysis and Design}, 47:453--460, 2011.

\bibitem{singhatprakash2011}
MK~Singha, T~Prakash, and M~Ganapathi.
\newblock Finite element analysis of functionally graded plates under
  transverse load.
\newblock {\em Finite Elements in Analysis and Design}, 47:453--460, 2011.

\bibitem{Sundararajan2005}
N~Sundararajan, T~Prakash, and M~Ganapathi.
\newblock Nonlinear free flexural vibrations of functionally graded rectangular
  and skew plates under thermal environments.
\newblock {\em Finite Elements in Analysis and Design}, 42(2):152--168, 2005.

\bibitem{valizadehnatarajan2013}
N~Valizadeh, S~Natarajan, OA~Gonzalez-Estrada, T~Rabczuk, TQ~Bui, and S~Bordas.
\newblock {NURBS-based finite element analysis of functionally graded plates:
  elastic bending, vibraton, buckling and flutter}.
\newblock {\em Composite Structures}, 99:209--326, 2013.

\bibitem{wu2004}
L~Wu.
\newblock Thermal buckling of a simply supported moderately thick rectangular
  {FGM} plate.
\newblock {\em Composite Structures}, 64:211--218, 2004.

\bibitem{wuliu2010}
SC~Wu, GR~Liu, XY~Cui, TT~Ngyuen, and GY~Zhang.
\newblock An edge-based smoothed point interpolation method {ES-PIM} for heat
  transfer analysis of rapid manufacturing system.
\newblock {\em International Journal of Heat and Mass Transfer}, 53:1938--1950,
  2010.

\bibitem{Zenkour2006}
AM~Zenkour.
\newblock Generalized shear deformation theory for bending analysis of
  functionally graded plates.
\newblock {\em Applied Mathematical Modeling}, 30:67--84, 2006.

\bibitem{Zenkour2007}
AM~Zenkour.
\newblock Benchmark trigonometric and 3{D} elasticity solutions for an
  exponentially graded thick rectangular plate.
\newblock {\em Archive of Applied Mechanics}, 77:197--214, 2007.

\bibitem{zenkourmashat2010}
AM~Zenkour and DS~Mashat.
\newblock Thermal buckling analysis of ceramic-metal functionally graded
  plates.
\newblock {\em Natural Science.}, 2:968--978, 2010.

\bibitem{zhaolee2009}
X~Zhao, YY~Lee, and KM~Liew.
\newblock Free vibration analysis of functionally graded plates using the
  element free $kp-${Ritz} method.
\newblock {\em Journal of Sound and Vibration}, 319:918--939, 2009.

\bibitem{zhaolee2009a}
X~Zhao, YY~Lee, and KM~Liew.
\newblock Mechanical and thermal buckling analysis of functionally graded
  plates.
\newblock {\em Composite Structures}, 90:161--171, 2009.

\end{thebibliography}

\end{document}